\documentclass{amsart}
\usepackage{amssymb, amscd, amsmath, amsthm, epsf, epsfig, latexsym}
 
\newtheorem{theorem}{Theorem}
\newtheorem{lemma}[theorem]{Lemma}
\newtheorem{corollary}[theorem]{Corollary}

\theoremstyle{definition}

\numberwithin{equation}{section}

\begin{document}

\title[Concordance Crosscap Numbers]{Concordance Crosscap Numbers of Knots and the Alexander Polynomial}
\author{Charles Livingston}

\maketitle

\begin{abstract}  For a knot $K$ the concordance crosscap number, $c(K)$, is the minimum crosscap number among all knots concordant to $K$.  Building on work of G. Zhang  which studied the determinants of knots with $c(K)  < 2$, we apply the Alexander polynomial to construct new algebraic obstructions to $c(K) < 2$.  With the exception of low crossing number knots previously known to have $c(K) < 2$, the obstruction applies to all but four prime knots of 11 or fewer crossings.
\end{abstract}
\vskip.3in

Every knot $K \subset S^3$ bounds an embedded surface $F \subset S^3$ with $F \cong \#_nP^2 - B^2$  for some $n \ge 0$, where $P^2$ denotes the real projective plane .  The crosscap number of $K$, $\gamma(K)$, is defined to be the  minimum such $n$.   The careful study of this invariant began with the work of Clark in~\cite{cl}; other references include~\cite{ht, my1}.    The study of the 4--dimensional crosscap number, $\gamma_4(K)$, defined   similarly but in terms of $F \subset B^4$,  appears  in such articles as~\cite{my2, vi, ya}.  

 Gengyu Zhang~\cite{zh} recently introduced a new knot invariant, the {\it concordance crosscap number}, $\gamma_c(K)$.  This is defined to be  the minimum crosscap number of any knot concordant to $K$.   This invariant is the nonorientable version of the concordance genus, originally studied by Nakanishi~\cite{na} and Casson~\cite{ca}, and later investigated in~\cite{li}.

In~\cite{zh}, Zhang presented an obstruction to $\gamma_c(K) \le 1$ based on the homology of the 2--fold branched cover of the knot, or equivalently, det$(K)$.  Inspired by her work, in this note we will observe that the obstruction found in~\cite{zh}  extends to one based on the Alexander polynomial of $K$, $\Delta_K(t)$, and the signature of $K$, $\sigma(K)$.

\begin{theorem}\label{mainthm} Suppose $\gamma_c(K) = 1$ and set $q = |\sigma(K)| +1$.  For all odd prime power divisors   $p$ of $q$, the $2p$--cyclotomic polynomial $\phi_{2p}(t)$ has odd exponent in $\Delta_K(t)$.  Furthermore, every other symmetric irreducible polynomial $\delta(t)$   with odd exponent in $\Delta_K(t)$ satisfies $\delta(-1) = \pm 1$. 

\end{theorem}

\begin{proof}   Any knot $K'$ with $\gamma(K') =1$ bounds a Mobius band, and is thus  a $(2,r)$--cable of some knot $J$ for some odd $r$.  If $K$ is concordant to $K'$, then $\sigma(K) = \sigma(K') = \pm(|r| -1)$; the signature $\sigma(K')$ is given by a formula of Shinohara~\cite{sh} for the signature of 2--stranded cable knots.  It follows that $|\sigma(K)| = |r| -1$, so $|r| = |\sigma(K) |+1 = q$.

According the a result of Seifert~\cite{se}, the Alexander polynomial 
of $K'$ is given by $\Delta_{2,q}(t)\Delta_J(t^2)$, where
 $\Delta_{2,q}(t)$ is the Alexander polynomial of the $(2,q)$--torus knot.  
A standard result states that  $\Delta_{2,q}(t) = \frac{(t^{2q}-1)(t-1)}{ ( t^{q} -1)(t^2-1)} = \frac{t^{q} +1}{t+1}$.  
This can be written as the product of cyclotomic polynomials, $$\Delta_{2,q}(t)= \prod_{p | q\ ,\ p >1} \Phi_{2p}(t).$$

Since $K$ is concordant to $K'$, $K \# -K'$ is slice, and thus has Alexander polynomial of the form $g(t)g(t^{-1})$.  That is, with $q = |\sigma(K)| +1$, 

$$\Delta_K(t) \Delta_J(t^2) \Delta_{2,q}(t) = g(t)g(t^{-1}). $$

 We now make two observations: {\bf (1)} Any symmetric irreducible polynomial has even exponent in $g(t)g(t^{-1})$, and thus even exponent in $ \Delta_{K}(t) \Delta_{J}(t^2) \Delta_{(2,q)}(t)$; {\bf (2) } since $\Delta_J(t)$ is an Alexander polynomial, $\Delta_J(1) = \pm 1$, 
 and thus $\Delta_J(t^2)|_{_{t=-1}} = \pm 1$.
 
 By Lemma~\ref{cyclemma}, $\Phi_{2p}(-1) = p$ if $p$ is an odd prime power, and  $\Phi_{2p}(-1) = \pm 1$ if $p$ is an  odd composite.   Thus, for $p$ an odd prime power divisor of $q$,  $\Phi_{2p}(t)$ has odd exponent in  $\Delta_{2,q}(t)$ and does not divide  $\Delta_{J}(t^2)$, so has odd exponent in  $\Delta_{K}(t)$.  Any other irreducible factor of $\Delta_{K}(t)$ with odd exponent is either a factor $\delta(t)$ of $\Delta_{2,q}(t)$, and thus of the form  $\Phi_{2p}(t)$ with $p$ an  odd composite (and so $\delta(-1) = \pm1$), or else is not a factor of $\Delta_{2,q}(t)$ and so has odd exponent in $\Delta_{J}(t^2)$, and again must satisfy $\delta(-1) = \pm 1$.  This completes the argument.
\end{proof}

\begin{lemma}\label{cyclemma} The cyclotomic polynomial $\Phi_{2p}(t)$ satisfies $ \Phi_{2p}(-1) = p$ if $p$ is an odd prime power and $ \Phi_{2p}(-1) = \pm1$ if $p$ is an odd composite.
\end{lemma}

\begin{proof} For an odd $r$, $h_r(t) = \frac{t^r + 1}{t +1}$ satisfies $h_r(-1) = r$ by L'Hospital's rule. We have that $h_r(t)$ is the product $$ h_r(t) = \prod_{p | r\ ,\ p >1} \Phi_{2p}(t).$$  For $p$ a prime power, $s^n$,    $\Phi_{2p}(t) = \frac{t^{s^n} +1}{t^{s^{n-1}} +1}$, and so, again by L'Hospital's rule, $\Phi_{2s^{n}}(-1) = s$.  Thus,   the product $$ \prod_{p | r\ ,\ p >1,\ p  \text{\ a prime power}} \Phi_{2p}(-1) = r. $$ It follows that all the other terms in the product expansion of $h_r(t)$ must equal $\pm 1$ when evaluated at $t=-1$, as desired.
\end{proof}

\vskip.2in 

\noindent{\bf Examples}  

Theorem~\ref{mainthm} is quite effective in ruling out $\gamma_c(K) = 1$.  For instance, there are 801 prime knots with 11 or fewer crossings.  Of these, 51 are known to be topologically slice, and 23 are known to be concordant to a $(2,q)$-torus knot for some $q$ and thus have $\gamma_c = 1$.  Of the remaining 727 knots, all but four can be shown to have $\gamma_c \ge 2$.  These four are $11n_{45}$ and $11n_{145}$, both of which are possibly slice, and $9_{40}$ and $11n_{66}$, both of which are possible concordant to the trefoil.  
Of the collection of 727 knots, Yasuhara's result~\cite{ya} applies to show that 207 of them have 4--ball crosscap number  $\gamma_4(K) \ge 2$.  The 4--ball crosscap numbers of the rest are unknown.

As a second set of examples, consider knots $K$ with $\Delta_K(t)$ of degree 2.  It follows immediately from Theorem~\ref{mainthm} that there are only two possibilities:  either $\sigma(K) = 0$ and $\Delta_K$ is reducible (an irreducible symmetric quadratic $f(t)$ cannot satisfy $f(1) = \pm1$ and $f(-1) = \pm 1)$, or $\sigma(K) = \pm 2$ and $\Delta_K(t) = t^2 -t+1$.  

We conclude with the further special case consisting of  $(p,q,r)$--pretzel knots, $P(p,q,r)$, with $p, q$, and $r$ odd; some of these were studied in~\cite{zh}.  If we let $D = D(p,q,r) = pq + qr +rp$, then   $$\Delta_{P(p,q,r)}(t) = \frac{D+1}{4}t^2  - \frac{D-1}{2}t + \frac{D+1}{4},$$ which has discriminant $-D$.  Thus, by the previous argument   we have: 

\begin{corollary} If $\gamma_c(P(p,q,r)) = 1$ then either $\sigma(P(p,q,r)) = 0$ and $D(p,q,r) = -l^2$ for some   integer $l$ or $\sigma(P(p,q,r) )= \pm 2$  and $D(p,q,r) = 3$.\end{corollary}

These pretzel knots include some shown by Zhang~\cite{zh} to have 4-dimensional crosscap number $\gamma_4(K) = 1$.

\newcommand{\etalchar}[1]{$^{#1}$} 

\begin{thebibliography}{BVS{\etalchar{+}}92}


\bibitem[Ca]{ca} A. Casson, 
unpublished.

\bibitem[Cl]{cl} B. Clark, {\sl 
Crosscaps and knots,}
Internat. J. Math. Math. Sci. {\bf  1} (1978), no. 1, 113--123. 

\bibitem[HT]{ht} M. Hirasawa and M. Teragaito, {\sl  Crosscap numbers of 2-bridge knots},  Topology {\bf 45}  (2006),  no. 3, 513--530.

 \bibitem[Li]{li} C. Livingston, {\sl The concordance genus of knots}, Alg. and Geom. Top. {\bf 4} (2004), 1--22.
 
 
\bibitem[MY1]{my1} H. Murakami and A. Yasuhara, {\sl Crosscap number of a knot},
Pacific J. Math. {\bf 171} (1995), no. 1, 261--273. 


\bibitem[MY2]{my2} H. Murakami and A. Yasuhara, {\sl  Four-genus and four-dimensional clasp number of a knot},  Proc. Amer. Math. Soc.  {\bf 128}  (2000),  no. 12, 3693--3699.

\bibitem[Na]{na} Y. Nakanishi, 
{\it A note on unknotting number},
  Math. Sem. Notes Kobe Univ. {\bf 9} (1981),  99--108.
 
 \bibitem[Se]{se} H. Seifert, {\sl  
On the homology invariants of knots}, Quart. J. Math., Oxford Ser. (2) {\bf 1},  (1950), 23--32. 
 
\bibitem[Sh]{sh} Y. Shinohara, {\sl 
On the signature of knots and links},
Trans. Amer. Math. Soc. {\bf 156}  (1971), 273--285. 

\bibitem[Vi]{vi} O.~Viro, {\sl Positioning in codimension 2, and the boundary}, Uspehi Mat. Nauk {\bf 30} (1975), 231--232.

\bibitem[Ya]{ya} A. Yasuhara, {\sl Connecting lemmas and representing homology classes of simply connected $4$-manifolds}, Tokyo J.~Math. {\bf 19} (1996), no. 1, 245--261. 


\bibitem[Zh]{zh} G.~Zhang, {\sl Concordance Crosscap Number of a Knot}, Bulletin of the London Mathematical Society {\bf 39} , 755-761.


\end{thebibliography}
\end{document}